\documentclass[12pt]{article}

 \newcommand{\nc}{\newcommand}
 \nc{\R}{{\rm I\!R}}
 \nc{\N}{{\rm I\!N}}
 \nc{\Prob}{{\rm I\!P}}
 \nc{\E}{{\rm I\!E}}

\begin{document}

\begin{center}
{\Huge {\bf Perfect Simulation from the\\
\vspace{0.2cm}
Quicksort Limit Distribution\\}}

\vspace{0.8cm}
\noindent 
{ Luc Devroye\\
School of Computer Science\\
McGill University\\
{\tt luc@cs.mcgill.ca}\\
{\tt http://www-cgrl.cs.mcgill.ca/\~{}luc/}\\

\vspace{0.8cm}
James Allen Fill\\
Department of Mathematical Sciences\\
The Johns Hopkins University\\ 
{\tt jimfill@jhu.edu}\\
{\tt http://www.mts.jhu.edu/\~{}fill/}\\

\vspace{0.8cm}
Ralph Neininger\\
Institut f\"ur Mathematische Stochastik\\
Universit\"at Freiburg\\
{\tt rn@stochastik.uni-freiburg.de}\\
{\tt http://www.stochastik.uni-freiburg.de/homepages/neininger/}}

 \end{center}

\vspace{0.3cm}

\begin{center}

{\bf Abstract}

\end{center}

\begin{quote}

The weak limit of the normalized number of comparisons needed by the {\tt Quicksort} algorithm to sort $n$
randomly permuted items is known to be determined implicitly by a distributional fixed-point equation. We
give an algorithm for perfect random variate generation from this distribution. 
\end{quote}
\vspace{0.3cm}
\noindent
{\small {\it Key words and phrases:\/}\ Quicksort, random variate generation, simulation, perfect simulation,
rejection method, Monte Carlo method, fixed-point equation.\\

}\\
\pagebreak

\noindent
\section{Introduction}
Let $C_n$ denote the number of key comparisons needed to sort a list of~$n$ randomly permuted items
by {\tt Quicksort}.  It is known that
$$
\E C_n = 2 (n + 1) H_n -4 n \sim 2 n \ln n\mbox{\ \ and\ \ Var}\,C_n \sim (7 - (2 \pi^2 / 3)) n^2, 
$$
where $H_n$ denotes the
$n$th harmonic number.  Furthermore,
\begin{eqnarray*}
X_n := \frac{C_n - \E C_n}{n} \longrightarrow X
\end{eqnarray*}
in distribution.  This limit theorem was first obtained by R\'egnier~\cite{reg} by an application of the
martingale convergence theorem.  R\"osler~\cite{roe} gave a different proof of this limit law via the
contraction method.  R\"osler's approach identifies the distribution of $X$ to be the unique solution
with zero mean and finite variance of the distributional fixed-point equation
\begin{eqnarray}
\label{fpe}
X \stackrel{\cal D}{=} UX^{(1)} + (1 - U) X^{(2)} + g(U),
\end{eqnarray}
where $X^{(1)}$, $X^{(2)}$, and $U$ are independent; $X^{(1)}$ and $X^{(2)}$ are distributed as $X$; $U$ is
uniform $[0,1]$; $g$ is given by $g(u) := 1 + 2 u \ln u + 2 (1 - u) \ln(1 - u)$; and $\stackrel{\cal D}{=}$
denotes equality in distribution.

The limit random variable $X$ has finite moments of every order which are computable from the fixed point
equation~(\ref{fpe}).  Tan and Hadjicostas~\cite{taha} proved that $X$ has a Lebesgue density.  Not much else
was known rigorously about this distribution until Fill and Janson recently derived some
properties of the limiting density~\cite{fija1} and results about the rate of convergence of the law
of~$X_n$ to that of~$X$~\cite{fija2}.  Some of these results are restated for the reader's convenience in the
next section.

We develop an algorithm, based on the results of Fill and Janson, which returns a perfect sample
of the limit random variable $X$.  We assume that we have available an infinite sequence of i.i.d.\ 
uniform~$[0,1]$ random variables. Our solution is based on a modified rejection method, where we use a
convergent sequence of approximations for the density to decide the outcome of a rejection test. Such an
approach was recently used by Devroye~\cite{dev2} to sample perfectly from perpetuities. 

\section{Properties of the quicksort density}

Our rejection sampling algorithm is based on a simple upper bound and an approximation of (the unique
continuous version of) the {\tt Quicksort} limit density~$f$.  We use the following properties of~$f$
established in~\cite{fija1} and~\cite{fija2}.  Let~$F_n$ denote the distribution function for~$X_n$.
\medskip

\noindent
P1. $f$ is bounded~\cite{fija1}:
\begin{eqnarray*}
\sup_{x\in\R} f(x) \le K:= 16,
\end{eqnarray*}
P2. $f$ is infinitely differentiable and~\cite{fija1}
\begin{eqnarray*}
\sup_{x\in\R} |f^\prime(x)| \le \widetilde{K}:= 2466,
\end{eqnarray*}
P3. With $\delta_n := (2 \hat{c} / \widetilde{K})^{1/2} n^{-1/6}$, where $\hat{c}:=(54 c K^2)^{1 / 3}$,
$c := 589$, we have~\cite{fija2}
\begin{eqnarray*}
\sup_{x \in \R} \left| \frac{F_n(x + (\delta_n / 2)) - F_n(x - (\delta_n/2))}{\delta_n} - f(x) \right|
  \le R_n,
\end{eqnarray*}
where $R_n:=(432 c K^2 \widetilde{K}^3)^{1/6} n^{-1/6}$.\\

By property P2, $f$ is Lipschitz continuous with Lipschitz constant $\widetilde{K}$.  Therefore, Theorem~3.5
in Devroye \cite[p. 320]{dev1} implies the upper bound 
\begin{eqnarray*}
f(x) \le \sqrt{2 \widetilde{K} \min(F(x), 1 - F(x))}.
\end{eqnarray*}
Here, $F$ denotes the distribution function corresponding to $f$.
Markov's inequality yields
$F(x)= \Prob(X \le x)\le (\E X^4) / x^4$ for all $x<0$. Similarly, $1 - F(x) = \Prob(X > x) \le (\E X^4)
/ x^4$ for $x > 0$.  The fourth moment of $X$ can be derived explicitly in terms of the zeta function
either by Hennequin's formula for the cumulants of $X$ (this formula was conjectured in Hennequin~\cite{hen}
and proved later in his thesis) or through the fixed point equation~(\ref{fpe}). From (\ref{fpe}), Cramer
\cite{cra} computed $\E X^4 = 0.7379 \ldots$ (accurate to the indicated precision), so $\E X^4 < 1$.
Therefore, if we define
\begin{equation}
\label{gdef}
g(x):=\min \left( K, (2 \widetilde{K})^{1/2} x^{-2} \right), \qquad x\in\R,
\end{equation}
we have $f \le g$. The scaled version $\widetilde{g}:=\xi g$ is the density of a probability measure
for $\xi := 1/ \|g\|_{L^1} = [4 K^{1/2} (2\widetilde{K})^{1/4}]^{-1}$.  A perfect sample from the
density~$\widetilde{g}$ is given by $[(2 \widetilde{K})^{1/4} / K^{1/2}] S U_1 / U_2$, with $S$, $U_1$,
and $U_2$ independent; $U_1$ and $U_2$ uniform $[0,1]$; and $S$ an equiprobable random sign (cf.\ Theorem~3.3
in Devroye \cite[p. 315]{dev1}).
\medskip

{\em Remark.\/}\ According to the results of~\cite{fija1}, $f$ enjoys superpolynomial decay at $\pm \infty$,
so certainly $f \leq g$ for some~$g$ of the form $g(x) := \min(K, C x^{-2})$.  One way to obtain an explicit
constant~$C$ is to use
$$
x^2 f(x) \leq \frac{1}{2 \pi} \int^{\infty}_{- \infty} |\phi''(t)|\,dt, \qquad x \in \R,
$$
where~$\phi$ is the characteristic function corresponding to~$f$, and to bound $|\phi''(t)|$ [e.g.,\ by
$\min(c_1, c_2 t^{-2})$ for suitable constants $c_1, c_2$] as explained in the proof of Theorem~2.9
in~\cite{fija1}.  But we find that our approach is just as straightforward, and gives a smaller value of~$C$
(although we have made no attempt to find the best~$C$ possible using the Fourier techniques
of~\cite{fija1}).

\section{The rejection algorithm}

We have found an explicit, integrable upper bound on $f$.  Furthermore, an approximation of $f$ with
explicit error estimate is given by P3.  Let
\begin{eqnarray*}
f_n(x) := \frac{F_n(x + (\delta_n / 2)) - F_n(x - (\delta_n / 2))}{\delta_n}
\end{eqnarray*}
with $\delta_n$ given in P3. Then $|f_n(x)-f(x)|\le R_n$ for all $x\in\R$, and $R_n\rightarrow 0$ for
$n\to\infty$. 

To calculate the values of $f_n$ we require knowledge about the probabilities of the events $\{C_n = i\}$.
Let $N(n, i)$ denote the number of permutations of $n$ distinct numbers for which {\tt Quicksort} needs
exactly $i$ key comparisons to sort. Then
\begin{eqnarray*}
\Prob(C_n = i)=\frac{N(n, i)}{n!}.
\end{eqnarray*}
These probabilities are non-zero only if $n-1 \le i \le n (n-1) / 2$. With the initializing conventions
$N(0,0) := 1$ and $N(i,0) := 0$ for $i \ge 1$ and the obvious values $N(n, i)=0$ for $i < n-1$ and
for $i > n (n-1) / 2$, we have the following recursion for $n \ge 1$ and for $n-1 \le i\le n (n - 1) / 2$:
\begin{eqnarray*}
N(n, i) = \sum_{k = 1}^n \sum_{l = 0}^{i - (n - 1)} N(k - 1, l) N(n - k, i -(n - 1) - l).
\end{eqnarray*}
This recurrence is well known.
To verify it, assume that the first pivot element is the $k$th largest element out of $n$. Then
the number of permutations leading to $i$ key comparisons is the number $N(k - 1,l)$ of permutations of the
items less than the pivot element which are sorted with $l$ key comparisons, multiplied by the corresponding
number of permutations for the elements greater than the pivot element, summed over all possible values of
$k$ and~$l$. Note that $n - 1$ key comparisons are used for the splitting procedure. Observe that we also
have $\E C_n = \sum_i i N(n, i) / n!$. The table $(N(n, i): i \le n (n - 1) / 2)$, and $\E C_n$, can be
computed from the previous tables $(N(k, i): i \le k (k - 1) / 2)$, $0 \leq k < n$, in time $O(n^5)$. Then,
observe that, for $y < z$,
\begin{eqnarray*}
F_n(z) - F_n(y) = \frac{1}{n!} \sum_{\E C_n + n y < i \le \E C_n + n z} N(n, i),
\end{eqnarray*}
and thus $f_n(x)$ is computable from the table $(N(n, i): i \le n (n - 1) / 2)$ and $\E C_n$ in time $O(n (z
- y)) = O(n \delta_n)=O(n^{5 / 6})$. Now, the following rejection algorithm gives a perfect sample~$X$ from
the {\tt Quicksort} limit distribution~$F$:\\

\noindent
{\tt repeat\\
\phantom{aa}generate $U$, $U_1$, $U_2$ uniform $[0,1]$\\
\phantom{aa}generate $S$ uniform on $\{-1, +1\}$\\
\phantom{aa}$X \leftarrow ((2\widetilde{K})^{1/4}/K^{1/2})SU_1/U_2$\\
\phantom{aa}$T \leftarrow U g(X)\;\;\;\;\;$(where $g(x) := \min(K,(2\widetilde{K})^{1/2}/x^2)$)\\
\phantom{aa}$n \leftarrow 0$\\
\phantom{aa}repeat\\
\phantom{aaaa}$n \leftarrow n+1$\\
\phantom{aaaa}compute the full table of $N(n, i)$ for all $i\le n(n-1)/2$\\
\phantom{aaaa}$Y \leftarrow f_n(X)$\\
\phantom{aa}until $|T-Y|\ge R_n$\\
\phantom{aa}Accept = [$T\le Y-R_n$]\\
until Accept\\
return $X$}\\

This algorithm halts with probability one, and produces a perfect sample from the {\tt Quicksort} limit
distribution.  The expected number of outer loops is $\|g\|_{L^1}=4K^{1/2}(2\widetilde{K})^{1/4} \doteq
134.1$.  Note, however, that the constants $K$ and~$\widetilde{K}$ are very crude upper bounds for
$\|f\|_\infty$ and $\|f^\prime\|_\infty$, which from the results of numerical calculations reported
in~\cite{taha} appear to be on the order of~$1$ and~$2$, respectively. 

Moreover, considerable speed-up could be achieved for our algorithm by finding another approximation $f_n$ to
$f$ that either is faster to compute or is faster to converge to $f$ (or both).  One promising approach, on
which we hope to report more fully in future work, is to let $f_1, f_2, \ldots$ be the densities one obtains,
starting from a suitably nice density~$f_0$ (say, standard normal), by applying the method of successive
substitutions to~(\ref{fpe}).  Indeed, Fill and Janson~\cite{fija2} show that then $f_n \to f$ uniformly at
an exponential rate.  However, one difficulty is that these computations require repeated numerical
integration, but it should be possible to bound the errors in the numerical integrations using calculations
similar to those in~\cite{fija1}.

\medskip

{\em Remark.\/}\ Let $k \equiv k_n := \lfloor \log_2(n + 1) \rfloor$.  We noted above that if $N(n, i) > 0$,
then $n - 1 \leq i \leq n (n - 1) / 2$.  This observation can be refined.  In fact, using arguments as
in~\cite{fil}, it can be shown that $N(n, i) > 0$ if and only if $m_n \leq i \leq M_n$, with
\begin{eqnarray*}
m_n &:=& k (n + 1) - 2^{k + 1} + 2 \sim n \log_2 n = (1 / \ln 2)\,n \ln n = (1.44\ldots)\,n \ln n\\
    & =& \mbox{the total path length for the complete tree on~$n$ nodes}
\end{eqnarray*}
and $M_n := n (n - 1) / 2$.  These extreme values satisfy the initial conditions $m_0 = 0 = M_0$ and, for
$n \geq 1$, the simple recurrences
$$
m_n = m_{n - 1} + \lfloor \log_2 n \rfloor \qquad \mbox{and} \qquad M_n = M_{n - 1} + (n - 1).
$$

\end{document}